# NONSTATIONARY RELAXED MULTISPLITTING METHODS FOR SOLVING LINEAR COMPLEMENTARITY PROBLEMS WITH H−MATRICES

Cuiyu Liu and Chenliang Li

School of Mathematics and Computing Science, Guangxi Colleges and Universities Key Laboratory of Data Analysis and Computation, Guilin University of Electronic Technology, Guilin, Guangxi, China,541004.

## *ABSTRACT*

*In this paper we consider some non stationary relaxed synchronous and asynchronous multisplitting methods for solving the linear complementarity problems with their coefficient matrices being H−matrices. The convergence theorems of the methods are given,and the efficiency is shown by numerical tests.*

## 1. INTRODUCTION

Many science and engineering problems are usually induced as linear complementarity problems(LCP): finding an $x \in R^n$ such that

$$x \geq 0, Ax - f \geq 0, x^{\bullet} (Ax - f) = 0, \quad (1.1)$$

where $A \in R^{n \times n}$ is a given matrix, and $f \in R^n$ is a vector. It is necessary to establish an efficient algorithm for solving the complementarity problem(CP). There have been lots of works on the solution of the linear complementarity problem([9,10,14,15,13,18]), which presented feasible and essential techniques for LCP.

The multisplitting method was introduced by O'Leary and White [17] and further studied by many people [11,12,1,2,3,4,5,6]. In the standard multisplitting method each local approximation solution $x^{k+1}$ is updated once using the same vector $x^k$. At the $k$ th iteration of a nonstationary multi splitting method, each processor $i$ solves the problem $q(k,i)$ times, in each time using the new obtained vector to update the $x^k$. [16] presented the following non-stationary multi splitting algorithm for linear systems:

**Algorithm 1.**(Nonstationarymultisplitting). Given the initial vector $x^0$,

For $k = 0,1,...$ until convergence
In processor $i$, $i = 1$ to $m$

$$y_i^0 = x^k,$$

For $l = 1$ to $q(k,i)$

$$F_i y_i^l = G_i y_i^{l-1} + b$$

$$x^{k+1} = \sum_{i=1}^{m} E_i y_i^{q(k,i)}.$$

In [16], relaxed nonstationarymultisplitting methods are also studied. The computational results show that these method are better than the standard multisplitting methods. [8] presented a nonstationary two-stage multisplitting methods with overlapping blocks. [7] proved the convergence of the nonstationarymultisplitting method for solving a system of linear equations when the coefficient matrix is symmetric positive definite.
The purpose of this paper is also on establishing efficient parallel iterative methods for solving the LCP. By skillfully using the matrix multisplitting methodology and the block property, we propose a class of nonstationarymultisplitting methods, for solving the linear complementarity problems (1.1).
The paper is organized as follows. In Section 2 we propose synchronous nonstationarymultisplitting method for solving LCP and establish its convergence theorem. In Section 3 we give an asynchronous nonstationary parallel multisplitting method for solving LCP and analysis the convergence of the algorithm. In Section 4, some numerical results show the efficiency of our algorithms.

## 2. Synchronous Relaxed Nonstationary multisplitting Method

Machida([13]) extended the multisplitting methods to the symmetric LCP. And Bai ([1,2,3,4]) developed a class of synchronous relaxed multisplitting methods for LCP, in which the system matrix is an $H$-matrix. In this section, by using multisplitting and block property techniques, we present a nonstationarymultisplitting method for the LCP(1.1), in which $A$ is an $H-$matrix.
At first we briefly describe the notations. In $R^n$ and $R^{n\times n}$ the relation $\geq$ denotes the natural components partial ordering.In addition, for $x, y \in R^n$ we write $x > y$ if $x_i > y_i, i = 1,2,\ldots,n$. A nonsingular matrix $A = (a_{ij}) \in R^{n\times n}$ is termed $M-$matrix,if $a_{ij} \leq 0$ for $i \neq j$ and $A^{-1} \geq 0$. Its comparison matrix $< A >= (\alpha_{ij})$ is defined by $\alpha_{ii} = | a_{ii} |$, $\alpha_{ij} = -| a_{ij} | (i \neq j)$. $A$ is said to be an $H-$matrix if $< A >$ is an $M-$matrix. $A$ is said to be an $H_+-$matrix if $A$ is an $H-$matrix with positive diagonal elements.

**Definition 2.1.** $A$ splitting $A = M - N$ is termed $M-$splitting of matrix $A$ if $M$ is an $M-$matrix and $N \geq 0$.

**Remark 2.1.** Let's separate the index set $S := \{1, 2, ..., n\}$ into $m$ nonempty subsets $S_i (i = 1, 2, ..., m)$ such that $\cup_{i=1}^{m} S_i = S$. We define that

$$E_{l,i} = (e_{pq}) = \begin{cases} \beta_{l,i} > 0, & p = q \in S_i, \\ 0, & otherwise, \end{cases}$$

$$\text{where } \sum_{i=1}^{m} E_{l,i} = I, \text{ and } E_{l,i} > 0.$$

According to the block property, some variables corresponding to the zero entries of $E_{l,i}$ need not be calculated.

In this section, we'll discuss a synchronous nonstationarymultisplitting method.

**Algorithm 2.** (Synchronous relaxed nonstationarymultisplitting method)

1) Give an initial value $x^0$, and let $k = 0$.

2) For each $i$ ($i = 1, 2, \ldots, m$),

$$y^{0,i} = x^k.$$

For $j = 1$ to $s(k,i)$,

$$\begin{cases} y^{j,i} \geq 0, \\ M_i y^{j,i} \geq F^{j,i}, \\ \left(y^{j,i}\right)^{\mathrm{T}} \left(M_i y^{j,i} - F^{j,i}\right) = 0, \end{cases} \tag{2.1}$$

where $F^{j,i} = f + N_i y^{j-1,i}$.

3)

$$x^{k+1} = \omega \sum_{i=1}^{m} E_i y^{s(k,i),i} + (1-\omega) x^k, \tag{2.2}$$

$$\text{where } \sum_{i=1}^{m} E_i = I, \ E_i > 0.$$

4) $k := k + 1$, return to step 2) .

**Lemma 2.1.** ([11]) Let $A, B \in R^{n \times n}$ satisfy that $\langle A\rangle \leq \langle B\rangle$. If $A$ is an $H$-matrix, then $B$ is also an $H$-matrix.

Following from Lemma 2.1, we can get the following lemma.

**Lemma 2.2.** If $\langle A \rangle \leq \langle M \rangle - |N|$, and $A$ is an $H$-matrix, then $M$ is an $H$-matrix.

**Lemma 2.3.** [11] Let $A = (a_{ij}) \in R^{n \times n}$. If there exists a $u \in R^n$, $u > 0$, such that $|A| u < u$, then there exists a number $\theta \in [0,1)$, such that $\rho(A) \leq \theta$.

**Lemma 2.4.** If $\langle A \rangle \leq \langle M \rangle - |N|$, and $A$ is an $H$-matrix, then $\rho(\langle M \rangle^{-1} |N|) < 1$.

**Proof:** By Lemma 2.2, $\langle M \rangle$ is an $M-$ matrix. Therefore, there exists a positive vector

$u=\langle A\rangle^{-1}e$, $e=(1,1,\ldots,1)^{\bullet}\in R^n$, such that

$$\langle M\rangle^{-1}|N|u\le(I-\langle M\rangle^{-1}\langle A\rangle)u=u-\langle M\rangle^{-1}e<u.$$

By Lemma 2.3, $\rho(\langle M\rangle^{-1}|N|)<1$.

The next lemma is obvious.

**Lemma 2.5.** Let $A\in R^{n\times n}$ be an $H_+$-matrix, and $x\in R^n$. If $x_j\ge 0$, then

$$(\langle A\rangle|x|)_j\le(Ax)_j.$$

**Lemma 2.6.** Let $A$ be an $H_+$ matrix, $A\le M_i-N_i$ ( $i=1,2,\ldots,m$ ) satisfy that $\langle A\rangle\le\langle M_i\rangle-|N_i|$ for each $i$, and $x^*$ be the solution of problem(1.1). If $y^{s(k,i),i}$ is generated by Algorithm 2, then

$$|y^{s(k,i),i}-x^*|\le\langle M_i\rangle^{-1}|N_i||y^{s(k,i)-1,i}-x^*|\le(\langle M_i\rangle^{-1}|N_i|)^{s(k,i)}|x^k-x^*|.$$

**Proof:** By Lemma 2.1 and Lemma 2.2, $\langle M_i\rangle$ is an $M-$ matrix. Consider the following cases:

(1) $y_j^{s(k,i),i}>x_j^*=0$.

By (1.1) and (2.1), we have

$$\left(M_ix^*-N_ix^*-f\right)_j\ge 0, (2.3)$$

and

$$\left(M_iy^{s(k,i),i}-N_iy^{s(k,i)-1,i}-f\right)_j=0. (2.4)$$

Substracting (2.4) from (2.3), we have

$$\left(M_i(y^{s(k,i),i}-x^*)\right)_j\le\left(N_i(y^{s(k,i)-1,i}-x^*)\right)_j\le\left(|N_i||y^{s(k,i)-1,i}-x^*|\right)_j.$$

Otherwise, by Lemma 2.5,

$$\left(M_i(y^{s(k,i),i}-x^*)\right)_j\ge\left(\langle M_i\rangle|y^{s(k,i),i}-x^*|\right)_j.$$

Therefore,

$$\left(\langle M_i\rangle|y^{s(k,i),i}-x^*|\right)_j\le\left(|N_i||y^{s(k,i)-1,i}-x^*|\right)_j. (2.5)$$

(2) $x_j^*>y_j^{s(k,i),i}=0$.

Similar to case (1), (2.5) holds true。

(3) $x_j^*>0$, $y_j^{s(k,i),i}>0$.

By (1.1) and (2.1), we have

$$\left(M_ix^*-N_ix^*-f\right)_j=0, (2.6)$$

and

$$\left(M_iy^{k,i}-N_iy^{s(k,i)-1,i}-f\right)_j=0. (2.7)$$

Substracting (2.7) from (2.6), we have

$$\left(M_i(y^{s(k,i),i}-x^*)\right)_j=\left(N_i(y^{s(k,i)-1,i}-x^*)\right)_j\le\left(|N_i||y^{s(k,i)-1,i}-x^*|\right)_j,$$

and

$$\left(M_i(y^{s(k,i),i}-x^*)\right)_j \ge \left(\langle M_i\rangle \mid y^{s(k,i),i}-x^* \mid\right)_j.$$

Therefore, (2.5) holds。

(4) If $y_j^{s(k,i),i}=x_j^*$, then $\left(x^*-y^{s(k,i),i}\right)_j=0$. From this we can deduce that the left side of (2.5) is non-positive, but the right side of (2.5) is non-negative. So (2.5) is true.

In a word,

$$\langle M_i\rangle \mid x^*-y^{s(k,i),i} \mid \le \mid N_i \mid\mid y^{s(k,i)-1,i}-x^* \mid,$$

or,

$$\mid x^*-y^{s(k,i),i} \mid \le \langle M_i\rangle^{-1} \mid N_i \mid\mid y^{s(k,i)-1,i}-x^* \mid.$$

Moreover, by induction,

$$\mid y^{s(k,i),i}-x^* \mid \le \langle M_i\rangle^{-1} \mid N_i \mid\mid y^{s(k,i)-1,i}-x^* \mid \le (\langle M_i\rangle^{-1} \mid N_i \mid)^{s(k,i)} \mid x^k-x^* \mid.$$

**Lemma 2.7.** ([11]) Let $A\in R^{n\times n}$ be nonsingular with $A^{-1}\ge 0$. Let $A=M-N=P-Q$ be two regular splittings of $A$ and

$$P^{-1}\ge M^{-1}.$$

Then

$$\rho(P^{-1}Q)\le \rho(M^{-1}N).$$

**Lemma 2.8.** Let $A=M-N$ and $A=D-B$ be $M-$ splittings of A, and $D=diag\{a_{11},a_{22},\cdots,a_{nn}\}$. If $M\le D$, then $\rho(M^{-1}N)\le\rho(D^{-1}B)<1$.

**Theorem 2.1.**Let $\gamma=\rho(D^{-1}B)$, $\omega\in(0,2/(1+\gamma))$. And let $A$ be an $H_+$-matrix, and for each $i=1,2,...,m$, $A=M_i-N_i$ satisfy that $\langle A\rangle\le\langle M_i\rangle-\mid N_i\mid$. Suppose that $\eta=\gamma/(\sum_{i=1}^{m}\left\|E_i\right\|)$, and $\tilde{s}$ satisfies that

$$\left\|(\langle M_i\rangle^{-1} \mid N_i \mid)^{s(k,i)}\right\|\le\eta,\quad s(k,i)\ge\tilde{s}, i=1,2,...,m. \tag{2.8}$$

If for each $k=1,2,...$, $i=1,2,\ldots,m$, $s(k,i)\ge\tilde{s}$, then the sequence $\{x^k\}$ generated by the Algorithm 2 converges to the solution $x^*$ of the problem (1.1).

**Proof:** Since $A$ is an $H_+$-matrix, and for each $i=1,2,\ldots,m$, $\langle A\rangle\le\langle M_i\rangle-\mid N_i\mid$. By Lemma 2.4, $\rho(\langle M_i\rangle^{-1}\mid N_i\mid)<1$. Let $T_k=\sum_{i=1}^{m}E_i(\langle M_i\rangle^{-1}\mid N_i\mid)^{s(k,i)}$. By Lemma 2.6, we have

$$\begin{aligned}\mid x^{k+1}-x^* \mid &\le \omega\sum_{i=1}^{m}E_i \mid y^{s(k,i),i}-x^* \mid + \mid 1-\omega \mid\mid x^k-x^* \mid\\ &\le(\omega T_k+\mid 1-\omega\mid I)\mid x^k-x^*\mid\\ &\le(\omega T_k+\mid 1-\omega\mid I)(\omega T_{k-1}+\mid 1-\omega\mid I)\cdots(\omega T_1+\mid 1-\omega\mid I)\mid x^0-x^*\mid.\end{aligned}$$

By (2.8), we have

$$\|T_k\| \le \sum_{i=1}^{m} \|E_i\| \|(\langle M_i \rangle^{-1} | N_i |)^{s(k,i)}\|$$
$$\le \gamma.$$

Therefore,

$$(\omega T_k + |1-\omega| I)|x^k - x^*| \le (\omega\gamma + |1-\omega|)|x^k - x^*|.$$

Moreover,

$$|x^{k+1} - x^*| \le (\omega T_k + |1-\omega| I)(\omega T_{k-1} + |1-\omega| I)\cdots(\omega T_1 + |1-\omega| I)|x^0 - x^*|$$
$$\le (\omega\gamma + |1-\omega|)^k |x^0 - x^*|.$$

As $0 < \omega < 2/(1+\gamma)$, then $\omega\gamma + |1-\omega| < 1$. Therefore, when $k \to \infty$, $|x^{k+1} - x^*| \to 0$.

Suppose that $A$ is an $M-$ matrix, it is an $H_+ -$ matrix, too. Therefore , we have the following corollary.

**Corollary 2.1.**Let $\theta = \rho(D^{-1}B)$, $\omega \in (0, 2/(1+\theta))$. And let $A$ be an $M$ -matrix, and for each $i = 1,2,...,m$, $A = M_i - N_i$ is an $M$ -splitting. Suppose that $\eta = \theta / (\sum_{i=1}^{m} \|E_i\|)$, and $\tilde{s}$ satisfies that

$$\|(M_i^{-1} N_i)^{s(k,i)}\| \le \eta, \quad s(k,i) \ge \tilde{s}, i = 1,2,...,m. \quad (2.9)$$

If for each $k = 1,2,...$, $i = 1,2,...,m$, $s(k,i) \ge \tilde{s}$, then the sequence $\{x^k\}$ generated by the Algorithm 2 converges to the solution $x^*$ of the problem (1.1).

# 3.ASYNCHRONOUS RELAXED NONSTATIONARY MULTISPLITTING METHOD

Bai([5,6]) proposed a class of standard asynchronous parallel multisplitting relaxation methods for LCP. Frommer ([12]) proposed an asynchronous weighted additive Schwarz scheme for solving the system of linear equations with multi-splitting method and Schwarz method. Mas [16] presented a relaxed nonstationarymultisplitting method for linear systems. Numerical experiments demonstrate the asynchronous method is faster than the corresponding synchronous one. The purpose of this section aims at extending this asynchronous nonstationary version to solving the LCP.

Let $N_0 = \{0,1,2,\ldots\}$. And for arbitrary $k \in N_0$, let $J(k) \subseteq \{1,2,\cdots,m\}$ be a nonempty set. Denote $J = \{J(k)\}_{k \in N_0}$ , $S = (s_1(k), s_2(k), \cdots, s_m(k))_{k \in N_0}$ as unbounded sequences, and have following properties:

(1) For any $i \in \{1,2,\cdots,m\}$, $k \in N_0$, $s_i(k) \le k.$

(2) For any $i \in \{1,2,\cdots,m\}$, $\lim_{k \to +\infty} s_i(k) = +\infty.$

(3) For any $i \in \{1,2,\cdots,m\}$, set $\{k \in N_0 \mid i \in J(k)\}$ is unbounded.

Let $s(k) = \min_i s_i(k), i = 1,2,\ldots,m$. It follows that $s(k) \le k$ by property (1). It's obvious that $\lim_{k \to +\infty} s(k) = +\infty$ by property (2).

Let $\{m_t\}_{t \in N_0}$ is a positive integer sequence which is strictly monotone increased.It is produced by following procedure:

$m_0$ is the least positive integer such that

$$\bigcup_{0 \le s(k) \le k < m_0} J(k) = \{1,2,\cdots,m\},$$

Similarly, $m_{t+1}$ is the least positive integer such that

$$\bigcup_{m_t \le s(k) \le k < m_{t+1}} J(k) = \{1,2,\cdots,m\}.$$

It's easily known that the sequence $\{m_t\}$ exists by the properties of the sets $J$ and $S$.

**Algorithm 3.** (Asynchronous relaxed nonstationarymultisplitting method)

1) Given an initial value $X^0 = \left(x^{0,1}, x^{0,2}, \ldots, x^{0,m}\right)^\bullet \in R^{nm}$, $x^{0,i} = x^0 \in R^n$, $i = 1,2,\ldots,m$. $k := 0$.

2) In processor $i$,

$$y^{i,0} = x^{s_i(k),i}.$$

For $v = 1$ to $q(i,k)$, $y^{i,v}$ is the solution of the following LCP:

$$\begin{cases} y^{i,v} \ge 0, \\ M_i y^{i,v} \ge F^{i,v}, \\ \left(y^{i,v}\right)^{\mathrm{T}}\left(M_i y^{i,v} - F^{i,v}\right) = 0, \end{cases} \qquad (3.1)$$

where $F^{i,v} = f + N_i y^{i,v-1}$.

Let

$$x^{k+1,l} = \begin{cases} x^{k,l}, & l \notin J(k), \\ \omega \sum_{i=1}^{m} E_{i,l}^{(k)} y^{i,q(i,k)} + (1-\omega) x^{k,l}, & l \in J(k), \end{cases}$$

where $\sum_{i=1}^{m} E_{i,l}^{(k)} = I$, $E_{i,l}^{(k)} > 0$, and $J(k) \subset \{1,2,\ldots,m\}$.

3) $k := k+1$, return to step 2).

The following lemma is obvious.

**Lemma 3.1.**Let $X^* = \left(x^*, x^*, \ldots, x^*\right)^\bullet$, $X^k = \left(x^{k,1}, x^{k,2}, \ldots, x^{k,m}\right)^\bullet$, if

$$\left\{X^{*}, X^{0}, X^{1}, \ldots, X^{k}\right\} \subset R^{mn}, \forall k \in N_{0},$$

and there exists a constant $\delta > 0$ and a positive vector $U = \left(u, u, \ldots, u\right)^{T} \in R^{mn}$, such that $| X^{q} - X^{*} | \leq \delta U$ for each $q \in \{0,1,2,\ldots,k\}$, then there holds

$$| v - X^{*} | \leq \delta U,$$

where

$$v = \begin{pmatrix} x^{s_1(k),1} \\ x^{s_2(k),2} \\ \vdots \\ x^{s_m(k),m} \end{pmatrix}$$

and $s_l\left(k\right) \leq k$ for all $l \in \{1,2,\ldots,m\}$.

**Theorem 3.1.**Let $\theta \leq \rho(D^{-1}B)$, $\omega \in (0, 2/(1+\theta))$. And let $A$ be an $H_{+}$-matrix, and for each $i = 1,2,\ldots,m$, $A = M_i - N_i$ satisfy that $\langle A\rangle \leq \langle M_i\rangle - | N_i |$. Suppose that $\eta = \theta / (\sum_{i=1}^{m} \left\| E_i \right\|)$, and $\tilde{q}$ satisfies that

$$\left\| (\langle M_i\rangle^{-1} | N_i |)^{q(i,k)} \right\| \leq \eta, \quad q(i,k) \geq \tilde{q}, i = 1,2,\ldots,m. \text{(3.2)}$$

If for each $k = 1,2,\ldots$, $i = 1,2,\ldots,m$, $q(i,k) \geq \tilde{q}$, then the sequence $\{x^{k,i}\}$ generated by Algorithm 3 converges to the solution $x^{*}$ of the problem (1.1).

**Proof:**Let $w = \langle A\rangle^{-1} e, e = \left(1,1,\ldots,1\right)^{\bullet} \in R^{n}$. Then, $w > 0$, and there exists a constant $\gamma \in \left[0,1\right)$, such that

$$\langle M_l\rangle^{-1} | N_l | w = (I - \langle M_l\rangle^{-1} \langle A\rangle) w = w - \langle M_l\rangle^{-1} e \leq \gamma w.$$

Let $w^{*} = \min_{1\leq j\leq n} \{w_j\}$. Then, $w^{*} > 0$. Define $\delta = \frac{\|x^{*} - x^{0}\|}{w^{*}}$, we get

$$| x^{*} - x^{0} | \leq \delta w.$$

Therefore, we have

$$\langle M_l\rangle^{-1} | N_l || x^{*} - x^{0} | \leq \delta \gamma w, \text{(3.3)}$$

and

$$| X^{0} - X^{*} | \leq \delta W,$$

where $W = (w, w, \ldots, w)^{\bullet}$.

Now we will prove that

$$| X^{k+1} - X^{*} | \leq \delta W, \quad k \in N_0. \qquad (3.4)$$

Let's assume that

$$| X^{q} - X^{*} | \leq \delta W, \quad q \in \{0,1,\ldots,k\}.$$

If $l \notin J\left(k\right)$, then $x^{k+1,l} = x^{k,l}$, and

$$|x^{k+1,l}-x^*|=|x^{k,l}-x^*|\le \delta w.$$

If $l\in J(k)$, then $x^{k+1,l}=\sum_{i=1}^{m}E_{i,l}^{(k)}y^{l,q(i,k)}$. Therefore,

$$|x^{k+1,l}-x^*|=|\omega\sum_{i=1}^{m}E_{i,l}^{(k)}y^{i,q(i,k)}+(1-\omega)x^{k,l}-x^*|$$

$$\le \omega\sum_{i=1}^{m}E_{i,l}^{(k)}|y^{i,q(i,k)}-x^*|+|1-\omega||x^{k,l}-x^*|. \quad (3.5)$$

By Lemma 2.6,

$$|y^{i,q(i,k)}-x^*|\le (\langle M_i\rangle^{-1}|N_i|)^{q(i,k)}|x^*-x^{s_i(k),i}|.$$

By Lemma 3.1, we have

$$|x^*-x^{s_i(k),i}|\le \delta w,$$

and then

$$|y^{i,q(i,k)}-x^*|\le (\langle M_i\rangle^{-1}|N_i|)^{q(i,k)}|x^*-x^{s_i(k),i}|. \quad (3.6)$$

Together (3.5) with (3.6), (3.2), and $\omega\theta+|1-\omega|<1$, we have

$$|x^{k+1,l}-x^*|\le \delta w.$$

This illustrates that for all $k\in N_0, |X^{k+1}-X^*|\le \delta W$.

In the sequel, we prove that

$$|X^k-X^*|\le \theta^t\delta W, \quad \forall k\ge m_t, t,k\in N_0. \quad (3.7)$$

By (3.4), we get

$$|X^k-X^*|\le \delta W=\theta^0\delta W.$$

Now we assume that , for any $k\ge m_t$, $|X^k-X^*|\le \theta^t\delta W$, and prove that (3.7) holds for any $k\ge m_{t+1}$.

By the definition of $m_{t+1}$, for any $k\ge m_{t+1}$ and $i\in\{1,2,\ldots,m\}$, there exists a positive integer $j$, satisfying $m_t\le s(j)\le j<k$, such that

$$x^{k,l}=x^{j+1,l}, l\in J(j),$$

where $x^{j+1,l}$ is the solution of (3.1).

Since for any $l\in\{1,2,...,m\}$, we have that $s(j)\le s_i(j)$. So we have $m_t\le s_i(j)$, and then $|x^*-x^{s_i(j),i}|\le \theta^t\delta w$. Therefore,

$$|x^{k,l}-x^*|=|x^{j+1,l}-x^*|$$

$$\le \sum_{i=1}^{m}E_{i,l}^{(j)}|\langle M_i\rangle^{-1}|N_i|||x^*-x^{s_i(j),i}|$$

$$\le \theta^{t+1}\delta w.$$

That is,

$$|X^k-X^*|\le \theta^{t+1}\delta W$$

holds for any $k\ge m_{t+1}$.

Since $\theta \in [0,1)$, we immediately get that $lim_{x\to\infty} x^{k,l} = x^*$, and then the sequence $\{x^{k,l}\}_{k\in N_0}$ generated by Algorithm 3 converges to the solution $x^*$ of problem (1.1).

**Corollary 3.1.**Let $\theta \le \rho(D^{-1}B)$, $\omega \in (0, 2/(1+\theta))$. And let $A$ be an $M$ -matrix, and for each $i = 1, 2, ..., m$, $A = M_i - N_i$ is an $M$ -splitting. Suppose that $\eta = \theta / (\sum_{i=1}^{m} \|E_i\|)$, and $\tilde{q}$ satisfies that

$$\|(M_i^{-1}N_i)^{q(i,k)}\| \le \eta, \quad q(i,k) \ge \tilde{q}, i = 1, 2, ..., m. \tag{3.8}$$

If for each $i = 1, 2, ..., m$, $k = 1, 2, \ldots$, $q(i,k) \ge \tilde{q}$, then the sequence $\{x^k\}$ generated by the Algorithm 3 converges to the solution $x^*$ of the problem (1.1).

## 4. Numerical Tests

In this section, we give some numerical results to illustrate the performance of the method presented in the paper. These results are for the purpose of illustrating new method for solving $LCP$ discussed in this paper. As we know, in practice the coefficient matrix is often a sparse matrix, so in the testing, we consider the $LCP$ as follows:

$$x \ge 0, Ax - b \ge 0, x^T(Ax - b) = 0,$$

where $b = (sin(2\pi/n), sin(4\pi/n), \cdots, sin(2\pi))^T$ is an $n \times 1$ vector,

$$A = \begin{pmatrix} C & -I & 0 & \cdots & 0 & 0 \\ -I & C & -I & \cdots & 0 & 0 \\ 0 & -I & C & \cdots & 0 & 0 \\ \vdots & \vdots & \vdots & \ddots & \vdots & \vdots \\ 0 & 0 & 0 & \cdots & C & -I \\ 0 & 0 & 0 & \cdots & -I & C \end{pmatrix}, C = \begin{pmatrix} 4 & -1 & 0 & \cdots & 0 & 0 \\ -1 & 4 & -1 & \cdots & 0 & 0 \\ 0 & -1 & 4 & \cdots & 0 & 0 \\ \vdots & \vdots & \vdots & \ddots & \vdots & \vdots \\ 0 & 0 & 0 & \cdots & 4 & -1 \\ 0 & 0 & 0 & \cdots & -1 & 4 \end{pmatrix},$$

and $I$ is a unit matrix.

For Case $m = 2$ and Case $m = 3$, we consider the nonstationary relaxed synchronous multisplitting method to solve the above linear complementarity problem. The stopping criterion of the out iteration is $\|x^{k+1} - x^k\| < 10^{-6}$. Let *time* denote the CPU time(sec.), and *Out-iter* denote the number of out iteration. The numerical results are listed in Table 1 and Table 2.

Table 1: $m = 2$, Comparison of results between Nonstationary Relaxed Multisplitting Method(NRM) and Standard Multisplitting Method(SMM)

| **level** | | **NRM** | | | | **SMM** |
|---|---|---|---|---|---|---|
| | | $\theta = 0.5$ | $\theta = 0.1$ | $\theta = 0.01$ | $\theta = 0.001$ | |
| **64*64** | **time** | **3.5100** | **4.4928** | **5.6472** | **7.5036** | **13.6657** |
| | **Out-iter** | **4378** | **3565** | **3409** | **3366** | **3360** |
| **128*128** | **time** | **44.0703** | **58.44** | **84.37** | **112.15** | **634.69** |
| | **Out-iter** | **17516** | **14262** | **13635** | **13462** | **13438** |
| **256*256** | **time** | **613.61** | **886.83** | **1439.0** | **2015.4** | **9492.8** |
| | **Out-iter** | **70082** | **57059** | **54549** | **53860** | **53761** |

Table 2: $m = 3$, Comparison of results between Nonstationary Relaxed Multisplitting Method(NRM) and Standard Multisplitting Method(SMM)

| level | | NRM | | | | SMM |
|---|---|---|---|---|---|---|
| | | $\theta = 0.5$ | $\theta = 0.1$ | $\theta = 0.01$ | $\theta = 0.001$ | |
| **64*64** | | 1.8876 | 2.4180 | 5.2416 | 9.7969 | 11.7157 |
| | Out-iter | 2055 | 1852 | 2721 | 4223 | 3360 |
| 128*128 | time | 19.7809 | 28.0022 | 51.7923 | 90.6210 | 166.0631 |
| | Out-iter | 7735 | 6639 | 7690 | 10814 | 13438 |
| **6*256** | time | 264.2813 | 386.3989 | 675.7807 | 978.6099 | 2767.2 |
| | Out-iter | 30365 | 25306 | 26410 | 30435 | 53761 |

From Table 1 and Table 2, the new nonstationary relaxed multisplittingmethod(NRM) has more efficiency than the standard multisplitting method(SMM). In our numerical tests, the stopping criterion of the inner iteration in SMM is $\left|\left(y^{i,v}\right)^{\mathrm{T}}\left(M_i y^{i,v} - F^{i,v}\right)\right| < 10^{-8}$. Though NRM has more number of the out iteration than SMM, the number of inner iteration of NRM is less than SMM. Therefore, NRM spend less time than SMM. These show us that NRM is more efficient.

## 5. Conclusion And Remarks

Mas ([16]) proposed a nonstationary parallel relaxed multisplitting methods for linear system. The numerical experiments show that these methods are better than the standard ones. In this paper, we develop a class of nonstationary relaxed synchronous and asynchronous multisplitting methods for solving linear complementarity problems with $H-$matrices. The convergence of the methods are analysed, and the efficiency is shown by numerical tests.

## ACKNOWLEDGEMENT

This work was supported by Natural Science Foundation of China (11161014),National Project for Research and Development of Major ScientificInstruments(61627807), and Guangxi Natural Science Foundation(2015 GXNSFAA 139014).